\documentclass{article}

\usepackage[nonatbib,final]{neurips_2023_ml4ps}

\usepackage[style=numeric,sorting=none]{biblatex}
\usepackage[utf8]{inputenc} 
\usepackage[T1]{fontenc}    
\usepackage{hyperref}       
\usepackage{url}            
\usepackage{booktabs}       
\usepackage{amsfonts}       
\usepackage{nicefrac}       
\usepackage{microtype}      
\usepackage{xcolor}         
\usepackage{algorithm}
\usepackage{algorithmicx}
\usepackage{algpseudocode}

\algdef{SE}[DOWHILE]{Do}{doWhile}{\algorithmicdo}[1]{\algorithmicwhile\ #1}%

\usepackage[utf8]{inputenc}
\usepackage[T1]{fontenc}
\usepackage{subfiles}
\usepackage{amsmath,amssymb,bm}
\usepackage{amsthm}
\usepackage{hyperref}
\usepackage{times}
\usepackage{framed}
\usepackage{graphicx}
\usepackage{svg}
\usepackage{verbatim}	
\usepackage{breqn} 
\usepackage{listings}
\usepackage{bbm}
\usepackage{subcaption}
\usepackage{xcolor}

\usepackage{url}            
\usepackage{booktabs}       
\usepackage{amsfonts}       
\usepackage{nicefrac}       
\usepackage{microtype}      
\usepackage{xcolor}         

\newcommand{\ee}{\end{equation}}
\newcommand{\be}{\begin{equation}}
\newcommand{\ec}{\end{center}}
\newcommand{\bc}{\begin{center}}
\newcommand{\eea}{\end{eqnarray}}
\newcommand{\bea}{\begin{eqnarray}}
\newcommand{\bd}{\begin{description}}
\newcommand{\ed}{\end{description}}
\newcommand{\bi}{\begin{itemize}}
\newcommand{\ei}{\end{itemize}}

\newcommand{\bx}{\bs{x}}

\newcommand{\bs}{\boldsymbol}

\newcommand{\bt}{\bs{\theta}}

\newcommand{\E}{\mathbb{E}}
\newcommand\calC{\mathcal{C}}

\newcommand\calL{\mathcal{L}}

\newcommand{\refeq}[1]{Eq. (\ref{#1})} 

\newcommand{\reffig}[1]{Fig. (\ref{#1})}

\title{Multi-fidelity Constrained Optimization for Stochastic Black Box Simulators}

\author{%
  Atul Agrawal\thanks{equal contribution} \,, \, Phaedon-Stelios Koutsourelakis \\
Professorship of Data-Driven Materials Modelling\\
Technical University of Munich\\
Munich, Germany\\
\texttt{\{atul.agrawal,p.s.koutsourelakis\}@tum.de} \\
   \And
  Kislaya Ravi\footnotemark[1]\,,  Hans-Joachim Bungartz \\
  Chair of Scientific Computing \\
  Technical University of Munich, Germany \\
  \texttt{\{kislaya.ravi,bungartz\}@tum.de} \\
}

\begin{document}

\maketitle

\begin{abstract}
Constrained optimization of the parameters of a simulator plays a crucial role in a design process. These problems become challenging when the simulator is stochastic, computationally expensive, and the parameter space is high-dimensional. One can efficiently perform optimization only by utilizing the gradient with respect to  the parameters, but these gradients are unavailable in many legacy, black-box codes. We introduce the algorithm Scout-Nd (\underline{S}tochastic \underline{C}onstrained \underline{O}p\underline{t}imization for N dimensions) to tackle the issues mentioned earlier by efficiently estimating the gradient, reducing the noise of the gradient estimator, and applying multi-fidelity schemes to further reduce computational effort. We validate our approach on standard benchmarks, demonstrating its effectiveness in optimizing parameters highlighting better performance compared to existing methods.
\end{abstract}

\section{Introduction}

Physics-based simulators are used across fields of engineering and science to drive research \cite{cranmer2020frontier}, and more recently used to generate synthetic training data for machine learning related tasks \cite{ros2016synthia}.
A common challenge is finding optimum parameters given some objective  subject to some constraints. High-dimensional parameter space and stochastic objective function make the optimization non-trivial.

Gradient-based methods have been shown to work well when the derivative is available \cite{degrave2019differentiable,agrawal2023probabilistic,de2018end,lucor2022simple}. However, for optimization/inference tasks involving physics-based simulators, only black-box evaluations of the objective are often possible(e.g., legacy solvers). It is commonly called simulation based inference(SBI)/optimization  \cite{cranmer2020frontier,ranganath2014black}. In such cases, one resort to gradient-free optimization \cite{more2009benchmarking}, for example, genetic algorithms \cite{banzhaf1998genetic} or Bayesian Optimization and their extensions \cite{snoek2012practical,menhorn2017trust}. 
The gradient-free methods perform poorly on high-dimensional parametric spaces \cite{more2009benchmarking}.
More recently, stochastic gradient estimators \cite{mohamed2020monte} have been used to estimate gradients of black-box functions and, hence, perform gradient-based optimization \cite{pflug2012optimization,louppe_adversarial_2019,shirobokov2020black,ruiz2018learning}. However, they do not account for the constraints.

This work introduces a novel approach for constrained stochastic optimization involving stochastic black-box simulators with high-dimensional parametric dependency. We draw inspiration from \cite{bird_stochastic_2018,staines_variational_2012} to estimate the gradients, extended it to include constraints and employed multi-fidelity strategies to limit the number of function calls, as the cost of running the simulator can be high. 
We choose popular gradient-free constrained optimization methods like Constrained Bayesian Optimization (cBO)\cite{gardner2014bayesian} and COBYLA \cite{powell1994direct}
to compare our method on standard benchmark problems.

\section{Methodology}

%
\textbf{Problem statement}
We are given a scalar valued function $f(\bx,\bm b)$ and a set of constraints $\calC(\bx,\bm b) = \{\calC_1(\bx,\bm b), \ldots, \calC_I(\bx,\bm b)\}$, where $\bx\in \mathbb{R}^d$ are the deterministic parameters and $\bm b$ represents a random vector \cite{pflug2012optimization}. The uncertainty may be caused by a lack of knowledge about the parameters or the inherent noise in the system. The objective $f$ or the constraints $\mathcal{C}$ depend implicitly on the output of the black-box simulator. Our task is to minimize the function $f(\bx,\bm b)$ with respect to $\bx$ subject to the constraints $\calC(\bx,\bm b)$. Because of the stochastic nature of the problem, we will optimize the objective function with respect to a robustness measure \cite{ben1999robust,bertsimas2011theory}. In this work, we will only consider the expectation as a robustness measure, in which case, the problem can be stated as  follows:
\begin{align}\label{eq:general_optimization}
	\min_{\bx} \E_{\bm b}[f(\bx,\bs{b})], \quad \textrm{s.t}~~~ \E_{\bm b}[\calC_i(\bx,\bs{b})]\leq 0 ,~~ \forall i \in \{1,\ldots,I\} 
\end{align}
In addition to that, the gradient of the objective function and the constraint is unavailable. Hence, one cannot directly apply gradient-based optimization methods.
\subsection{Constraint augmentation}\label{sec:constraints_stochastic}
We cast the constrained optimization problem (\refeq{eq:general_optimization}) to an unconstrained one using penalty-based methods \cite{wang_stochastic_2003, nocedal1999numerical}. We define an augmented objective function $\mathcal{L}$ as follows: 
\begin{align}\label{eq:updated_objective}
	\mathcal{L}(\bx, \bs{b},\bm{\lambda}) =  f(\bx,\bs{b})+ \sum_{i=i}^{I}\lambda_i \max\left(\calC_i(\bx,\bs{b}),0\right)
\end{align}
where $\lambda_i>0$ is the penalty parameter for the $i^{th}$ constraint and the $\max{(\cdot,\cdot)}$ controls the magnitude of penalty applied.
Incorporating the augmented objective (\refeq{eq:updated_objective}) in the \refeq{eq:general_optimization}, one can arrive at the following optimization problem: 
\begin{align}\label{eq:updated_opt_prob}
    \min_{\bx} \E_{\bs{b}} [\mathcal{L}(\bx, \bs{b},\bm{\lambda}) ]
\end{align}
The expectation is approximated by Monte Carlo which induces noise and necessitates stochastic optimization methods.
We alleviate the dependence on the penalty parameter $\bm \lambda$ by using the sequential unconstrained minimization technique (SUMT) algorithm \cite{fiacco1990nonlinear} where one starts with a small penalty term and gradually increases its value.
\subsection{Estimation of derivative}\label{sec:non-diff_obj_cons}
The direct computation of derivatives of $\mathcal{L}$ with respect to the optimization variables $\bm x$ is not feasible because of the unavailability of the gradients of the objective function and the constraints.
One notes many active research threads across disciplines which are trying to tackle this bottleneck \cite{cranmer2020frontier, louppe_adversarial_2019, beaumont2002approximate,marjoram2003markov, agrawal2023probabilistic,shirobokov2020black, shirobokov2020black}. We draw inspiration from  the Variational Optimization \cite{bird_stochastic_2018,staines2013optimization,ruiz2018learning}, which constructs an upper bound of the objective function as shown below:
\begin{align}
	\min \int \mathcal{L}(\bx,\bs{b},\bm{\lambda}) p(\bm b) d\bm{b} \leq \int \mathcal{L}(\bx,\bs{b},\bm{\lambda}) p(\bm b) q(\bx \mid \bt)d\bm{b}d\bx={U(\bt)}
\end{align}
where $q(\boldsymbol{x}\mid \bs{\theta})$ is a  density over the design variables $\bx$ with parameters $\bm{\theta}$. If $\bx^*$ yields the minimum of the objective $\E_{\bm b}[\calL]$, then this can be achieved with a degenerate $q$ that collapses to a Dirac-delta, i.e. if $q(\boldsymbol{x}\mid \bs{\theta})=\delta(\bx-\bx^*)$. The inequality above would generally be strict for all other densities $q$ or parameters $\bt$.  
Hence, instead of minimizing $\E_{\bm b}[\calL]$ with respect to $\bm{x}$, we can minimize the upper bound $U(\bt)$ with respect to the distribution parameters $\bm{\theta}$. Under mild restrictions outlined by \cite{staines_variational_2012}, the bound $U(\boldsymbol{ \theta})$ is differential w.r.t $\bm{\theta}$. One can evaluate the gradient of $U(\bt)$ as shown below:
\begin{align}\label{eq:grad_estimator}
	\nabla_{\bt} U(\bt) 
	=\E_{\bx, \bm{b}}\left[\nabla_{\bt} \log q(\bx \mid \bt) \mathcal{L}(\bx,\bs{b},\bm{\lambda})\right]
\end{align}
The Monte Carlo estimation of the expectation shown in \refeq{eq:grad_estimator} is as follows:
\begin{align}\label{eq:gradient_estimate}
	\frac{\partial U}{\partial \theta} \approx \frac{1}{S} \sum_{i=1}^{S} \mathcal{L}(\bx_i,\bs{b}_i,\bm{\lambda}) \frac{\partial}{\partial \bt} \log q\left(\bx_i \mid \theta\right)
\end{align}
The \refeq{eq:gradient_estimate} is known as the score function estimator \cite{glynn1990likelihood}, which also appears in the context of reinforcement learning \cite{williams1992simple}. 
The gradient estimation can be computationally expensive as each step will involve calling the simulator $S$ number of times. This step is can be easily parallelized.
\subsection{Variance reduction} \label{sec:grad_noise_reduction}
To reduce the mean square error of the estimator in \refeq{eq:gradient_estimate}, we propose the use of baseline discussed in \cite{kool_buy_2022} as shown below:
\begin{align}
	\frac{\partial U}{\partial \theta} \approx \frac{1}{S} \sum_{i=1}^{S}  \frac{\partial}{\partial \bt} \log q\left(\bx_i \mid \bt\right) \left(\mathcal{L}(\bx_i,\bs{b}_i,\bm{\lambda}) - \frac{1}{S-1}\sum_{j=1, j \neq i}^S \mathcal{L}(\bx_j,\bs{b}_j,\bm{\lambda}) \right) \label{eq:baseline_trick}
\end{align}
The above is an unbiased estimator and implies no additional cost beyond the  $S$ samples.
We also propose to use Quasi-Monte Carlo (QMC) sampling \cite{dick2013high} for variance reduction. QMC replaces $S$ randomly drawn samples by a pseudo-random sequence of samples of length $S$ with low discrepancy. This sequence covers the underlying design space more evenly than the random samples, thereby reducing the variance of the gradient estimator. 
\reffig{fig:grad_var_reduction} numerically shows the advantage of using variance reduction technique. We observe that the variance of the gradient is decreased using the variance reduction methods, 
specifically in high dimensions where we  observe $\sim 10\times$ benefit.
\subsection{Multi-fidelity}
The main computational bottleneck of the gradient estimation using \refeq{eq:gradient_estimate} is the multiple evaluation of the objective function. This becomes a significant concern for computationally expensive simulators. We propose to solve this problem using the multi-fidelity (MF) method \cite{peherstorfer2018survey}. Suppose we are given a set of $L$ functions modeling the same quantity and arranged in ascending order of accuracy and computational cost $\{ f_1, f_2, \ldots, f_L \}$. We want to optimize the design parameter with respect to the highest-fidelity model ($f_L$). We can estimate the gradient of the corresponding upper bound using the method suggested in \cite{giles2015multilevel} as shown below:
\begin{align}\label{eq:mf_gradient_estimate}
	\frac{\partial U_L}{\partial \theta} \approx \sum_{\ell=1}^{L} \frac{1}{S_{\ell}} \sum_{i=1}^{S_{\ell}} \left( \mathcal{L}_{\ell}(\bx_i,\bs{b}_i,\bm{\lambda})  - \mathcal{L}_{\ell-1}(\bx_i,\bs{b}_i,\bm{\lambda}) \right) \frac{\partial}{\partial \bt} \log q\left(\bx_i \mid \theta\right)
\end{align}
where $S_{\ell}$ is the number of samples used in the estimator at level $\ell$ and $\mathcal{L}_0(\cdot) = 0$. We want to use more samples from the low-fidelity model and lesser samples as we increase the fidelity. The method to calculate the number of samples is discussed in \cite{giles2015multilevel}. 

\subsection{Implementation details}
In the present study, we use \texttt{PyTorch} \cite{paszke2019pytorch} to efficiently compute the gradient. After the gradient estimation, we use the ADAM optimizer \cite{kingma2014adam} as the stochastic gradient descent method. In this work, $q(\boldsymbol{x}\mid \bs{\theta})$ takes the form of a Gaussian distribution with parameters $\bs{\theta} = \{ \mu , \sigma\}$ representing mean and variance. 
Our proposed algorithms is summarized in Algorithm \ref{Alg:optimization_algo}. The source code will be made available upon publication.
\begin{algorithm}
	\caption{Stochastic constrained optimization for non-differentiable objective (Scout-Nd)}
	\label{Alg:optimization_algo}
	\begin{algorithmic}[1]
		\State \textbf{Input}: Objective function(s), constraint(s), distribution $q(\bm{x}\mid\bt)$, stochastic gradient descent optimizer $\mathcal{G}$ and its hyper-parameters $\eta$, list of penalty terms $\{\bm{\lambda}_1, \bm{\lambda}_2, \ldots, \bm{\lambda}_K \}$, $\bm{\lambda}_K \rightarrow \infty$
		\State $\bt_0^0 \gets$ choose starting point, $k \gets 1$
		\Do
		\State $n \gets 0$
		\Do
		\State $\bm{x}_i \sim q(\bm{x} \mid \bt_k^n), \quad \bm{b}_i \sim p(\bm b)$ \Comment{Sampling step} 
        \State Evaluate augmented objectives $\mathcal{L}(\bx_i, \bs{b}_i,\bm{\lambda}_k)$ \Comment{\refeq{eq:updated_objective}}
		\State Monte-Carlo gradient estimate $\nabla_{\bt} U$ \Comment{\refeq{eq:baseline_trick}, \ref{eq:mf_gradient_estimate}}
		\State $\bt_k^{n+1} \gets \mathcal{G}(\bt_k^n, \eta, \nabla_{\bt} U)$ and $n \gets n + 1$ \Comment{Stochastic Gradient Descent}
		\doWhile{$\|\bt_{k}^{n} - \bt_{k}^{n-1}\| > \varepsilon_{\theta}$} 
		\State  $\bt_{k+1}^{0} \gets \bt_{k}^{n}$; $\{\mu, \sigma\} \gets \bm{\theta}_k^n$ and $k \gets k+1$
		\doWhile{$\| \sigma \| > \varepsilon_{\sigma}$} \Comment{Collapse to Dirac-delta}
		\State \Return{$\{\mu, \sigma\}$}
	\end{algorithmic}
\end{algorithm}

\section{Numerical Illustrations}
\vspace{-0.5em}

\begin{figure*}[!htbp]
\begin{minipage}{\textwidth}
  \begin{minipage}[b]{0.32\textwidth}
        \centering
    	\includegraphics[width=0.9\textwidth]{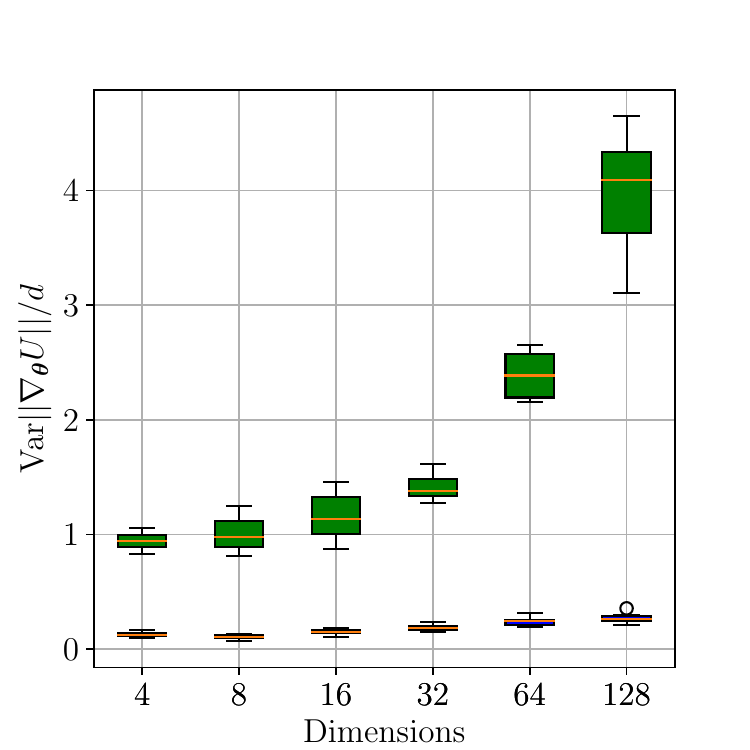}
    	\captionof{figure}{Box plot of the variance of the gradient estimator w.r.t the dimensions with 10 repeated runs for \refeq{eq:sphere_problem}. Gradient estimated with 128 samples at $\bt = \{\bs{1}_d,e^{\bs{1}_d}\}$. Green : no variance reduction, red : variance reduction}
    	\label{fig:grad_var_reduction}
    \end{minipage}
    \hfill
    \begin{minipage}[b]{0.32\textwidth} 
        \centering
    	\includegraphics[width=1\textwidth]{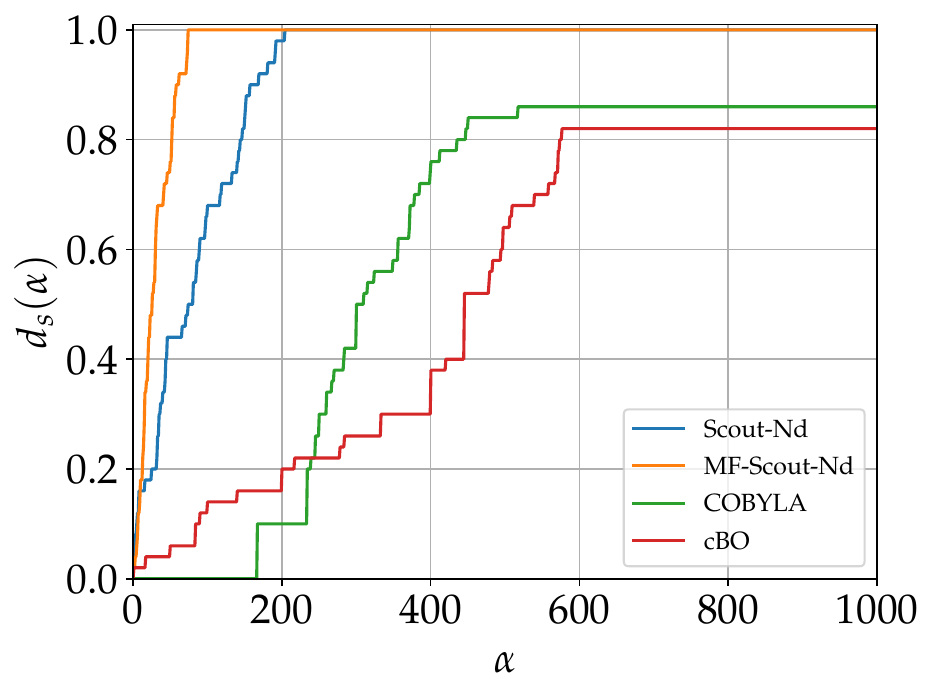}
    	\captionof{figure}{Data profiles plot for the set of optimizers. For the expectation estimation for each step, MF-Scout-Nd uses $S_1=10$ HF and $S_2 =50$ LF evaluations. Every other method uses 50 black-box evaluations. $\epsilon_f = 0.1$}
    	\label{fig:comparisions}
    \end{minipage}
    \hfill
    \begin{minipage}[b]{0.32\textwidth} 
        \centering
    \includegraphics[width=1\textwidth]{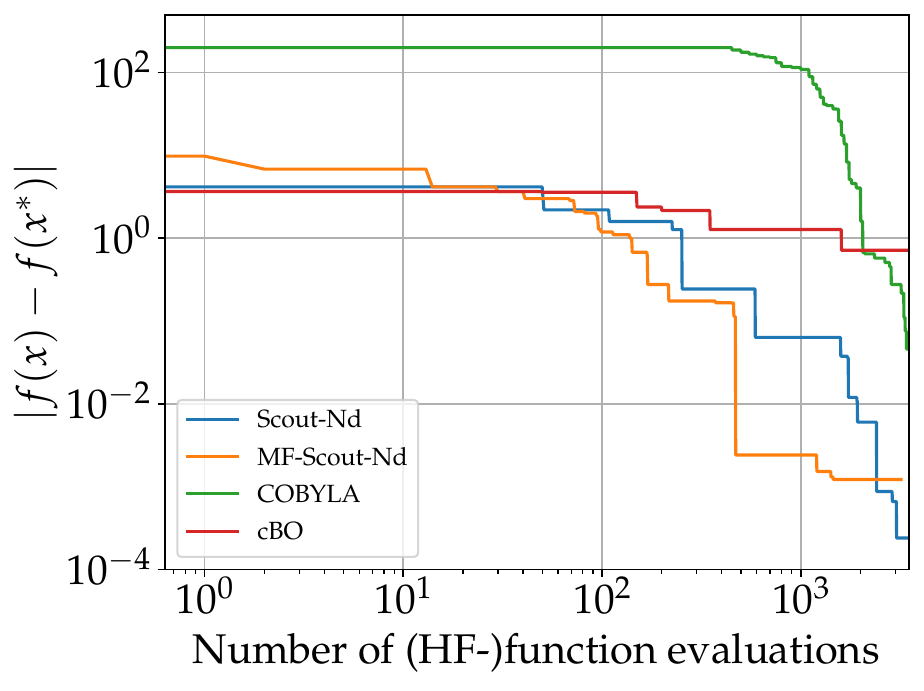}
    \captionof{figure}{Evolution of the distance between the optimal objective values suggested by the optimizer $(f(x))$ and the theoretical optimum $(f(x^*))$ with respect to the number of (HF-)function evaluations for the first case of sphere problem \ref{eq:sphere_problem} with $d=8$.}
    \label{fig:evolution_of_error}
    \end{minipage}
\end{minipage}
\end{figure*}

We discuss the numerical results of the proposed (MF)Scout-Nd algorithm on sphere-problem of varying dimensions ($d=\{2,4,8,16,32\}$). We use the data profiles proposed in \cite{more2009benchmarking} to compare (MF)Scout-Nd with cBO \cite{gardner2014bayesian} and COBYLA \cite{powell1994direct}, which are standard derivative-free optimization methods that can handle constraints.
We consider the following optimization problem with \textit{noisy} objective:
\begin{align}\label{eq:sphere_problem}
    \min_{\bx} \E_{b}\left[ \sum_{i=1}^{d} x_i^2 + b \right]; \quad s.t. \quad \mathcal{C}(\bm x)\leq 0
\end{align}
with $b\sim \mathcal{N}(0,0.1)$. We consider two different constraint cases. The first case is $\mathcal{C}(\bm x) = 1-(x_1+x_2)$ where the optimum lies on the constraint i.e. $\bm{x}^* = \{0.5, 0.5\} \cup \{0\}_{d-2}$. In the second case, the constraint surface is defined as $\mathcal{C}(\bm x) = \sum_{i=1}^{d} x_i-1$ leading to the optimum at $\bm{x}^* =\{0\}_{d}$.  
Let $\mathcal{S}$ be the set of optimizers $s$ and $\mathcal{P}$ be a set of benchmark problems $p$. Then the data profile $d_s(\alpha)$ \cite{more2009benchmarking} of a optimizer $s \in \mathcal{S}$  is given by 
\begin{align}\label{eq:data_profile}
    d_s(\alpha)=\frac{1}{|\mathcal{P}|} \left| \left\{p \in \mathcal{P}: \frac{t_{p, s}}{d_p+1} \leq \alpha\right\} \right|
\end{align}
where $d_p$ in the number of design variables in $p$, $\alpha$ is the allowed number of functional evaluations scaled by the number of design variables and $t_{p,s}$ is the minimum number of function calls a solver $s$ requires to reach the optimum of a problem $p$ within accuracy level $\epsilon_f$. We run each benchmark $5$ times leading to $|\mathcal{P}|= 5 (|d|) \times 5 (\text{number of runs}) \times 2 (\text{number of cases}) = 50$. 
For (MF)Scount-Nd, we consider two levels of multi-fidelity. The high-fidelity (HF) is given by the \refeq{eq:sphere_problem} and the low-fidelity (LF) is defined by down-scaling the $x_i$ in \refeq{eq:sphere_problem} to have $x_i = x_i/1.05$.

We can observe from \reffig{fig:comparisions} that Scout-Nd performs better than cBO and COBYLA because it solved most of the benchmark problems $p \in \mathcal{P}$ for a given $\alpha$. Our proposed algorithm outperforms the other two derivative-free methods because we use derivative approximation to move toward the optimum. This not only helped us to converge faster but also tackled high-dimensionality. MF-Scout-Nd performed better as it converged faster towards the optimum than Scout-Nd because it needs fewer costly function evaluations. We can also observe from \reffig{fig:evolution_of_error} that Scout-Nd and MF-Scout-Nd come closer to the actual optimal objective for a given computational budget.

\section{Conclusions}
\vspace{-0.5em}
We extended the method proposed by \cite{staines_variational_2012,bird_stochastic_2018} to account for constraints. We further employed a multi-fidelity strategy and gradient variance reduction schemes to reduce the number of costly simulator calls. We demonstrated on a classical benchmark problem that our method performs better than the chosen baselines in terms of quality of optimum and number of functional calls. As future work, we will test our method on real-world problems (for example: \cite{shirobokov2020black,louppe_adversarial_2019}).

\section{Broader Impact Statement}

Many real-world systems in engineering and physics are modeled by complex simulators that might be parameterized by a high-dimensional random variable. Some notable examples include particle physics, fluid mechanics, molecular dynamics, protein folding, cosmology, material sciences, etc. Frequently, the simulators are black-box, making the task of optimization/inference challenging, specifically in high dimensions. The task can be further complicated with the inclusion of constraints. In the present work, we introduced an algorithm to approximate the gradients for an optimization/inference task involving these simulators. We demonstrated that the proposed method performs better than the state-of-the-art on a standard benchmark problem. 

We do not see any direct ethical concerns associated with this research. The impact on society is
primarily through the over-arching context of providing novel methods to approach optimization/inference involving complex simulators.

\clearpage
\printbibliography

\end{document}